\documentclass[preprint,10pt]{elsarticle}

\usepackage[a4paper,top=4.0cm,bottom=4.0cm]{geometry}

\usepackage{amsmath,amssymb,amsthm}

\usepackage{graphicx}
\usepackage{subcaption}
\usepackage{tikz}
\usetikzlibrary{shapes,arrows.meta,calc,positioning,decorations.markings}

\usepackage{hyperref}
\hypersetup{colorlinks=true, urlcolor=red, citecolor=blue, linkcolor=red}

\makeatletter

\newcommand*\MyORCID{0009-0004-3236-9811}
\newcommand*\MyORCIDLink{%
  \href{https://orcid.org/\MyORCID}{\MyORCID}%
}

\gdef\emailauthor#1#2{%
  \stepcounter{ead}%
  \g@addto@macro\@elseads{%
    \raggedright
    \let\corref\@gobble
    \def\@@tmp{#1}%
    \eadsep{\ttfamily\expandafter\strip@prefix\meaning\@@tmp}%
    \space(#2, ORCID:\space\MyORCIDLink)%
    \def\eadsep{\unskip,\space}%
  }%
}

\makeatother

\begin{document}

\begin{frontmatter}

\title{A Surface-Based Formulation of the Traveling Salesman Problem}

\author{Y\i lmaz Arslano\u{g}lu}
\ead{yilmaz@hoketo.com}

\affiliation{organization={Independent Researcher}, country={Hamburg, Germany}}

\begin{abstract}
We present an exact formulation of the symmetric Traveling Salesman Problem (TSP) that replaces the classical edge-selection view with a surface-building approach.
Instead of selecting edges to form a cycle, the model selects a set of connected triangles where the boundary of the resulting surface forms the tour.
This method yields a mixed-integer linear programming (MILP) formulation where a tree constraint enforces global connectivity, while local connectivity at each vertex is guaranteed via Euler characteristic constraints, replacing the need for subtour elimination.
The formulation is exact when applied to the complete set of all triangles, despite being computationally intractable for all but the smallest instances. In practice, it provides a compact and effective heuristic when restricted to a sparse candidate set such as Delaunay triangulation.
\end{abstract}


\end{frontmatter}

\newcommand{\safeincludegraphics}[2][]{%
  \IfFileExists{#2}{\includegraphics[#1]{#2}}{\fbox{Missing image: #2}}%
}

\begin{figure}[h!]
\centering
\makebox[\textwidth][c]{\safeincludegraphics[width=1.10\textwidth, trim={2cm 2cm 2cm 2cm}, clip]{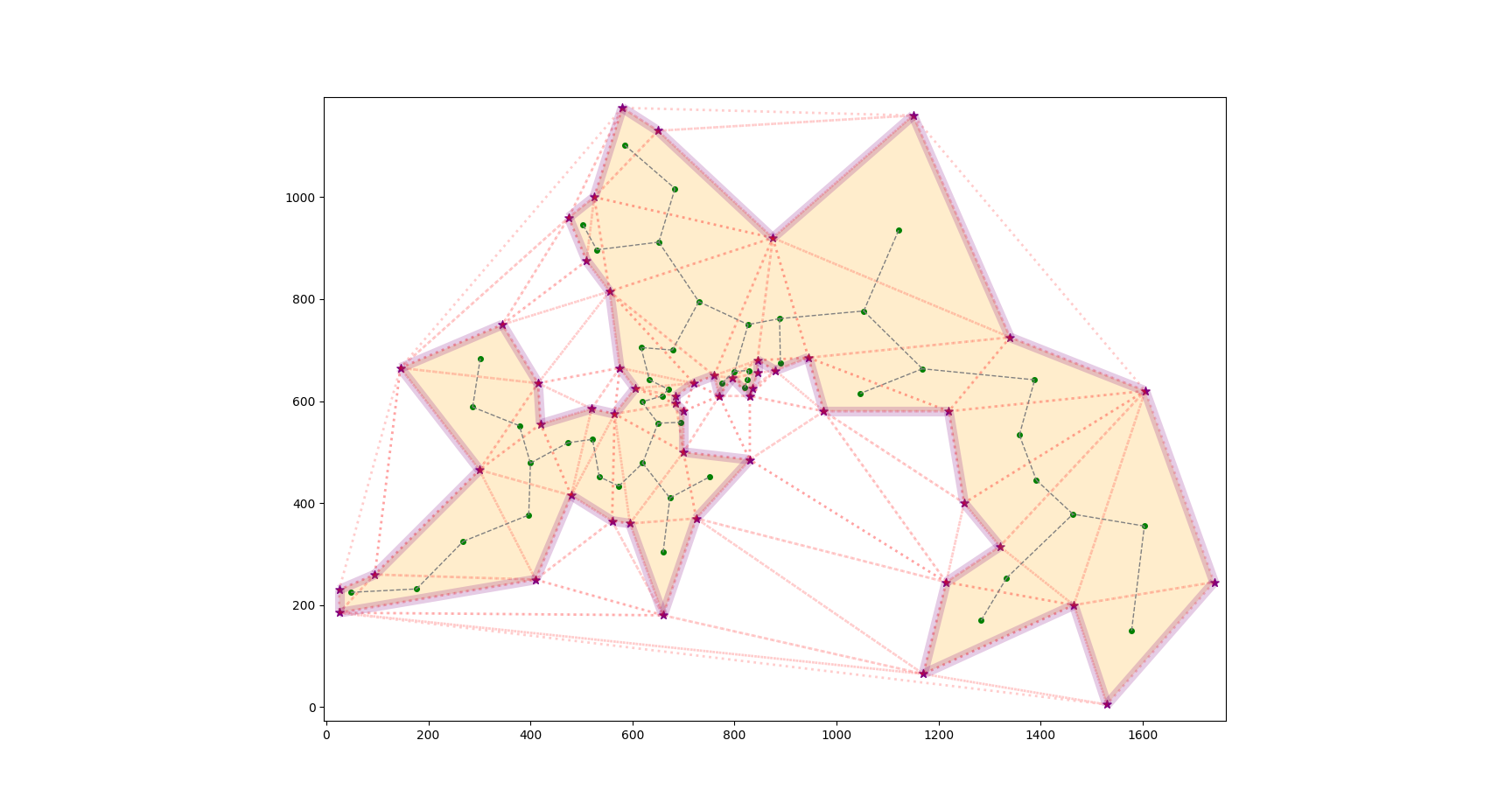}}
\caption{Visual intuition on \texttt{berlin52} ($N=52$): select a connected set of triangles (surface); its boundary is the optimal tour.}
\label{fig:teaser}
\end{figure}

\section{Introduction}
\begin{figure}[h!]
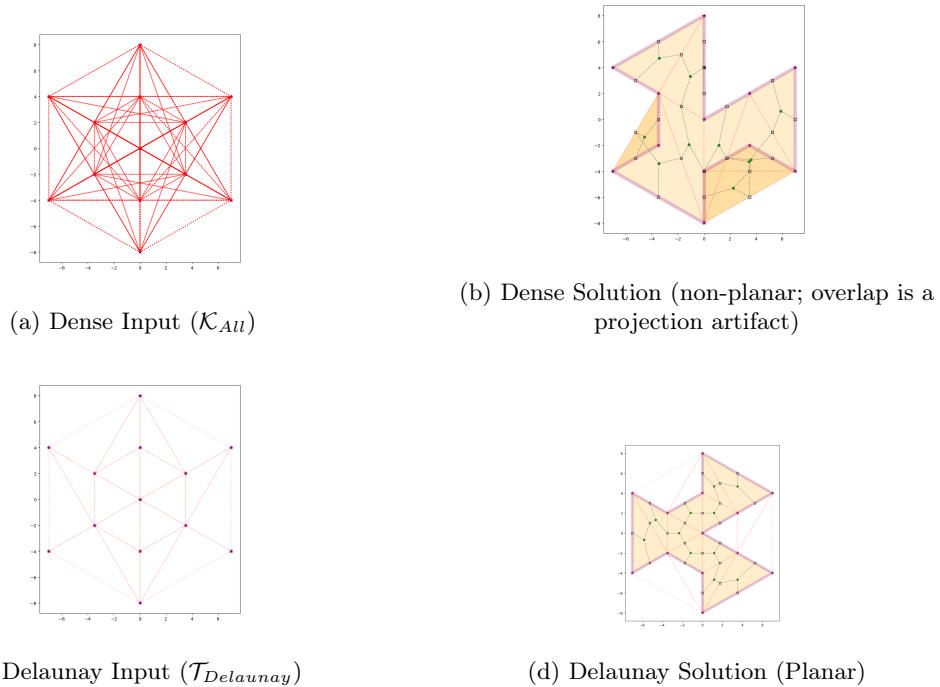

\centering \setlength{\tabcolsep}{1pt}
\begin{minipage}[b]{0.495\textwidth}
    \centering \safeincludegraphics[width=\textwidth]{nested\_hexagons\_all\_triangles\_input.png}
    \par\vspace{2pt} \small (a) Dense Input ($\mathcal{K}_{All}$)
\end{minipage}\hfill
\begin{minipage}[b]{0.495\textwidth}
    \centering \safeincludegraphics[width=\textwidth]{nested\_hexagons\_all\_triangles\_output.png}
    \par\vspace{2pt} \small (b) Dense Solution (non-planar; overlap is a projection artifact)
\end{minipage}\vspace{5pt}
\begin{minipage}[b]{0.495\textwidth}
    \centering \safeincludegraphics[width=\textwidth]{nested\_hexagons\_delaunay\_input.png}
    \par\vspace{2pt} \small (c) Delaunay Input ($\mathcal{T}_{Delaunay}$)
\end{minipage}\hfill
\begin{minipage}[b]{0.495\textwidth}
    \centering \safeincludegraphics[width=0.77\textwidth]{nested\_hexagons\_delaunay\_output.png}
    \par\vspace{2pt} \small (d) Delaunay Solution (Planar)
\end{minipage}
\caption{Nested hexagons ($N=13$): with all triangles the optimum is an abstract non-planar surface; Delaunay input forces a planar one.}
\label{fig:nested_hexagons}
\end{figure}

TSP is conventionally modeled on a graph $G=(V,E)$: choose a minimum-weight Hamiltonian cycle by selecting edges $E' \subseteq E$ (Fig.~\ref{fig:teaser}).
Compact edge-based models enforce local degree conditions but must also encode global connectivity, often via strong but exponential subtour elimination constraints~\cite{dantzig1954}, or via weaker compact alternatives such as MTZ~\cite{miller1960} or flow-based constraints~\cite{fox1980}.
We propose a shift in perspective from \textit{graph-theoretic} cycle selection to a \textit{topological} construction.
Instead of building the tour as a 1-dimensional sequence of edges, we construct a 2-dimensional surface whose boundary is the Hamiltonian cycle.
Specifically, given a candidate triangle set $\mathcal{T} \subseteq \binom{V}{3}$, we select a connected subset of triangles to minimize the total length of its induced boundary edges.
When $\mathcal{T}$ is the universal set $\mathcal{K}_{All}:=\binom{V}{3}$ (Fig.~\ref{fig:nested_hexagons}a), the formulation is exact but computationally feasible only for very small $N$.
Restricting $\mathcal{T}$ to a sparse triangulation (e.g., Delaunay - Fig.~\ref{fig:nested_hexagons}c) yields a highly efficient candidate-set heuristic, though it may render the model infeasible if $\mathcal{T}$ contains no Hamiltonian tour~\cite{dillencourt1990}.

\subsection{Scope}
Our goal is to provide a compact, exact alternative formulation that is structurally different from cycle-based models. We do not aim to compete with state-of-the-art branch-and-cut solvers for pure TSP, but rather to demonstrate that enforcing surface admissibility can yield strong relaxations and useful candidate-set models.

\subsection{Related Work}
Classical exact methods rely on branch-and-cut with subtour elimination~\cite{dantzig1954, applegate2006}.
Tree-structured connectivity constraints are well-studied in flow-based routing formulations~\cite{gavish1978traveling} and in Steiner variants~\cite{letchford2013}.
Geometric approximation schemes for Euclidean TSP exploit recursive planar subdivisions~\cite{arora1998, mitchell1999}.
Our formulation aligns with boundary-minimization viewpoints from computational topology~\cite{edelsbrunner2010} and geometry~\cite{fekete2016, fekete2022}, planar duality principles~\cite{whitney1931}, and
topological optimization of cycles within fixed surfaces~\cite{chambers2009}, but differs in that we optimize the \emph{surface itself} via a simplicial incidence structure.
This yields an incidence-graph formulation whose objective has a boundary-cancellation structure related in spirit to prize-collecting terms, except that the ``prizes'' arise from simplicial incidences, and is related in spirit to geometric covering variants such as TSP with Neighborhoods~\cite{dumitrescu2003}.

\subsection{Contribution}
\begin{figure}[t!]
\centering
\begin{tikzpicture}[scale=1.0, thick]
    \tikzset{
        city/.style={circle,fill=black,inner sep=1.5pt,outer sep=1pt},
        tri_node/.style={circle,draw=orange!80!black,fill=orange!30,inner sep=2pt},
        edge_node/.style={rectangle,draw=black,fill=white,inner sep=2.5pt},
        dual_edge/.style={draw=blue!80!black,dashed,line width=0.6pt},
        selected_dual/.style={draw=green!60!black,line width=1.6pt,opacity=0.55},
        tri_fill/.style={fill=orange!25,draw=none,opacity=0.28},
        boundary_edge/.style={draw=red!85!black,ultra thick}
    }

    \def\missingTri{1}

    \coordinate (C) at (0,0);
    \foreach \i in {1,...,6} \coordinate (V\i) at ({(\i-1)*60 + 90}:2);

    \foreach \i [evaluate=\i as \next using {int(mod(\i,6)+1)}] in {1,...,6} {
        \ifnum\i=\missingTri
        \else
            \path[tri_fill] (C) -- (V\i) -- (V\next) -- cycle;
        \fi
    }

    \draw[gray!30, line width=0.45pt]
        (C) -- (V1) (C) -- (V2) (C) -- (V3) (C) -- (V4) (C) -- (V5) (C) -- (V6);
    \draw[gray!30, line width=0.45pt]
        (V1) -- (V2) -- (V3) -- (V4) -- (V5) -- (V6) -- cycle;

    \foreach \i [evaluate=\i as \next using {int(mod(\i,6)+1)}] in {1,...,6} {
        \ifnum\i=\missingTri
        \else
            \draw[boundary_edge] (V\i) -- (V\next);
        \fi
    }
    \pgfmathtruncatemacro{\nextMissing}{int(mod(\missingTri,6)+1)}
    \draw[boundary_edge] (C) -- (V\missingTri);
    \draw[boundary_edge] (C) -- (V\nextMissing);

    \node[city, label={[font=\small, yshift=-2pt]below:$v$}] (c_center) at (C) {};
    \foreach \i in {1,...,6} \node[city] at (V\i) {};

    \foreach \i [evaluate=\i as \next using {int(mod(\i,6)+1)}] in {1,...,6}
        \node[tri_node] (t\i) at (barycentric cs:C=1,V\i=1,V\next=1) {};

    \foreach \i in {1,...,6}
        \node[edge_node] (e_in\i) at (barycentric cs:C=1,V\i=1) {};
    \foreach \i [evaluate=\i as \next using {int(mod(\i,6)+1)}] in {1,...,6}
        \node[edge_node] (e_out\i) at (barycentric cs:V\i=1,V\next=1) {};

    \foreach \i [evaluate=\i as \next using {int(mod(\i,6)+1)}] in {1,...,6} {
        \draw[dual_edge] (t\i) -- (e_in\i);
        \draw[dual_edge] (t\i) -- (e_in\next);
        \draw[dual_edge] (t\i) -- (e_out\i);
    }

    \foreach \i [evaluate=\i as \next using {int(mod(\i,6)+1)}] in {1,...,6} {
        \ifnum\i=\missingTri
        \else
            \draw[selected_dual] (t\i) -- (e_out\i);
            \draw[selected_dual] (t\i) -- (e_in\i);
            \draw[selected_dual] (t\i) -- (e_in\next);
        \fi
    }

    \begin{scope}[shift={(3.05, 0)}]

        \node[tri_node] at (0, 1.25) {};
        \node[right, font=\footnotesize, align=left] at (0.3, 1.25)
            {\textbf{Triangle node} ($t \in U$)\\(Cost: $0$)};

        \node[edge_node] at (0, 0.55) {};
        \node[right, font=\footnotesize, align=left] at (0.3, 0.55)
            {\textbf{Edge node} ($e \in W$)\\(Cost: $2L_e$)};

        \draw[dual_edge] (0, -0.10) -- (0.5, -0.10);
        \node[right, font=\footnotesize, align=left] at (0.6, -0.10)
            {Incidence edge ($\mathcal{A}$)\\(Profit: $L_e$)};

        \draw[selected_dual] (0, -0.55) -- (0.5, -0.55);
        \node[right, font=\footnotesize, align=left] at (0.6, -0.55)
            {Selected in tree};

        \draw[boundary_edge] (0, -1.00) -- (0.5, -1.00);
        \node[right, font=\footnotesize, align=left] at (0.6, -1.00)
            {Induced boundary};
    \end{scope}
\end{tikzpicture}
\caption{Bipartite incidence graph: triangles $x_t$ (circles) and edges $y_e$ (squares); dashed lines are adjacencies.
A feasible solution selects all but one candidate triangle (filled orange), inducing the highlighted boundary (red) around the center; translucent green edges indicate incidence edges chosen in the resulting tree.}
\label{fig:incidence_graph}
\end{figure}
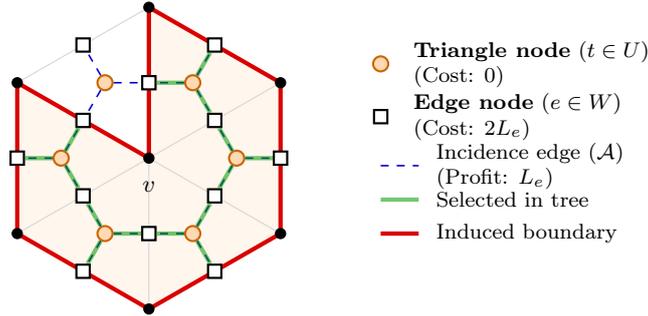

Our contribution is threefold:
\begin{enumerate}
    \renewcommand{\labelenumi}{(\roman{enumi})}
    \item We present an equivalence between min-weight Hamiltonian cycles in the \emph{complete} complex $\mathcal{K}_{All}$ (Fig.~\ref{fig:nested_hexagons}a) and max-weight dual trees in the bipartite incidence graph (Fig.~\ref{fig:incidence_graph});
    \item We give a MILP that replaces subtour elimination by a tree constraint plus an Euler characteristic filter (Sec.~\ref{sec:formulation});
    \item We demonstrate that the formulation is fully \emph{input-decoupled}: it treats the candidate triangle set $\mathcal{T}$ as a black-box constraint on the search space, working unchanged on the universal set $\mathcal{K}_{All}$ := $\binom{V}{3}$ as well as on sparse geometric triangulations (e.g., Delaunay or Greedy)~\cite{dillencourt1990, das1989, letchford2008}.
\end{enumerate}

Note: For geometric intuition, readers may first consider the planar case described in Appendix~\ref{apx:planar} (Fig.~\ref{fig:stokes_cancellation}) before diving into the general case.

\section{Topological Dual via the Incidence Graph}
This section introduces the combinatorial structure underlying the surface-based formulation.

\subsection{Incidence Graph}
Given a TSP instance on $G=(V,E)$, we construct the bipartite \emph{triangle--edge incidence graph}
$B=(U\cup W,\mathcal{A})$ (Fig.~\ref{fig:incidence_graph}), defined as follows:
\begin{itemize}
    \item \textbf{Triangle nodes ($U$):} a candidate triangle set $U := \mathcal{T} \subseteq \binom{V}{3}$.
    \item \textbf{Edge nodes ($W$):} the set of primal edges $W := E$.
    \item \textbf{Incidences ($\mathcal{A}$):} an arc $(t,e)\in\mathcal{A}$ exists if triangle $t\in U$ contains edge $e\in W$.
\end{itemize}
A feasible solution of our formulation corresponds to selecting a connected subset of triangles
together with the induced active edge set, represented as an active subgraph $B'\subseteq B$.

\subsection{Net Weight and Boundary Cancellation}
To model boundary minimization, we associate a cost with each primal edge node $e\in W$ and a
profit with each incidence $(t,e)\in\mathcal{A}$. For an active subgraph $B'$, define the net weight
\begin{equation}
W(B') := \sum_{(t,e)\in\mathcal{A}(B')} p(t,e) \;-\; \sum_{e\in V(B')\cap W} c(e). \label{eq:dual_score}
\end{equation}
With the choice $c(e)=2L_e$ and $p(t,e)=L_e$, internal edges (shared by two selected triangles)
cancel out, while boundary edges contribute exactly their length. As a result, optimizing $W(B')$
is equivalent to minimizing the total length of the induced boundary tour.~\cite{arslanoglu2026}.

\section{Mathematical Formulation}
\label{sec:formulation}
We formulate the problem as finding a subgraph in the bipartite incidence graph $B=(U \cup W, \mathcal{A})$ that corresponds to a valid surface. Let $x_t, y_e, z_{te} \in \{0,1\}$ be decision variables representing the selection of triangles, edges, and their incidences, respectively. The objective is to minimize the length of the boundary $\partial \mathcal{K}$ of the resulting surface complex $\mathcal{K}$:
\begin{equation}
    \text{Min } Z = \sum_{e \in W} 2 L_e y_e - \sum_{(t,e) \in \mathcal{A}} L_e z_{te} \label{eq:objective_main}
\end{equation}
Eq.~(\ref{eq:objective_main}) is exactly $Z=-W(B')$ for the induced active subgraph $B'$ (Eq.~\ref{eq:dual_score}), so minimizing $Z$ corresponds to maximizing the net weight.
The objective exhibits boundary cancellation and equals the exact length of the induced boundary tour. ~\cite{arslanoglu2026}.

\begin{figure}[t!]
\centering
\begin{subfigure}[b]{0.45\textwidth}
    \centering
    \begin{tikzpicture}[scale=0.9, line join=round]
        \node[anchor=south, font=\bfseries\footnotesize] at (0, 3.05) {A. Non-manifold edge};
        \coordinate (B) at (0, 0); \coordinate (T) at (0, 3);
        \filldraw[fill=orange!20, draw=orange!80!black, thin, opacity=0.8] (B)--(T)--(2, 1.5)--cycle;
        \filldraw[fill=orange!30, draw=orange!80!black, thin, opacity=0.8] (B)--(T)--(-1.5, 2)--cycle;
        \filldraw[fill=orange!50, draw=orange!80!black, thin, opacity=0.8] (B)--(T)--(-1.2, 0.5)--cycle;
        \draw[ultra thick, black] (B) -- (T);
        \node[rectangle, draw=black, fill=white, inner sep=2pt] at (0, 1.5) {};
    \end{tikzpicture}
\end{subfigure}
\hfill
\begin{subfigure}[b]{0.45\textwidth}
    \centering
    \begin{tikzpicture}[scale=0.7]
        \node[anchor=south, font=\bfseries\footnotesize] at (0, 3.05) {B. ``Bowtie'' singularity};

        \coordinate (C)  at (0,1.5);
        \coordinate (L1) at (-1.5,2.5);
        \coordinate (L2) at (-1.5,0.5);
        \coordinate (R1) at ( 1.5,2.5);
        \coordinate (R2) at ( 1.5,0.5);

        \fill[orange!20] (C) -- (L1) -- (L2) -- cycle;
        \fill[orange!20] (C) -- (R1) -- (R2) -- cycle;

        \fill[red] (C) circle (3pt);
        \node[right, red, font=\bfseries] at (0.1, 1.5) {$v$};

        \draw[magenta, thick] (C) -- (L1) -- (L2) -- cycle;
        \draw[magenta, thick] (C) -- (R2) -- (R1) -- cycle;
    \end{tikzpicture}
\end{subfigure}

\caption{Forbidden anomalies: (A) Edge incident to $>2$ triangles; (B) A "bowtie" vertex violating the Euler Filter.}
\label{fig:constraints_viz}
\end{figure}
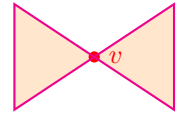

We enforce topological validity via four structural constraints (see Appendix \ref{apx:full_model} for the explicit MILP constraints):

\begin{enumerate}
    \item \textbf{Global Tree Connectivity:} $B'$ must form a single connected tree (cf. gray dashed lines connecting green nodes in Fig.~\ref{fig:teaser}). This ensures $\mathcal{K}$ is simply connected (topologically a disk).
    \item \textbf{Cardinality:} To guarantee $\mathcal{K}$ is a triangulation of a topological disk with $N$ vertices, we enforce $\sum x_t = N-2$ (triangles) and $\sum y_e = 2N-3$ (edges).
    \item \textbf{Manifold Regularity (Node Degrees):} Every primal edge must be incident to at most two triangles (see Fig.~\ref{fig:constraints_viz}A). In the dual, this restricts the degree of edge nodes $y_e$ within the incidence graph.
    \item \textbf{Local Connectivity (Euler Filter):} To prevent vertex singularities (e.g., ``bowties,'' Fig.~\ref{fig:constraints_viz}B), the neighborhood of each city $v$ must be connected (cf.\ the center vertex in Fig.~\ref{fig:incidence_graph}, where a boundary path arises only by omitting an incident triangle). We enforce this using the local Euler characteristic of the incidence subgraph $H_v\subseteq B'$~\cite{edelsbrunner2010}:
    \begin{equation}
    \chi(H_v)=|\mathrm{V}(H_v)|-|E(H_v)|=1 \quad \forall v\in V. \label{eq:euler_main}
    \end{equation}
    Here $H_v$ denotes the active subgraph of $B'$ induced by the triangles and edges incident to $v$ (see Appendix~\ref{apx:full_model}
    for the explicit definition). Under the global tree constraint, enforcing $\chi(H_v)=1$ leads to a simple path, ruling out
    disconnected vertex links and local ``bowtie'' singularities (Fig.~\ref{fig:constraints_viz}B)~\cite{arslanoglu2026}.

\end{enumerate}

\section{Theoretical Analysis}
\label{sec:theory}
All formal correctness proofs (boundary cancellation, disk admissibility, Hamiltonicity, and the resulting exactness of the MILP on the full complex) are given in the companion note~\cite{arslanoglu2026}.
Here we keep only the high-level takeaway: the objective measures the induced boundary length, and the tree+manifold+Euler constraints ensure the selected 2-complex is a disk whose boundary is a single Hamiltonian cycle.

\section{Computational Validation and Results}
We validated the correctness of the formulation on simple topological cases before measuring performance.

\subsection{Visual Verification (Topology)}
To provide a visual verification, we apply the model to an instance consisting of two nested hexagons with a central node ($N=13$). Fig.~\ref{fig:nested_hexagons} compares the solver's behavior under two different candidate triangle sets: the universal set $\mathcal{K}_{All}$ and the Delaunay triangulation $\mathcal{T}_{Delaunay}$.
In both cases, the formulation successfully extracts the optimal surface within the chosen input complex.

\subsection{Validation on Non-Metric Spaces}
\begin{figure}[htbp]
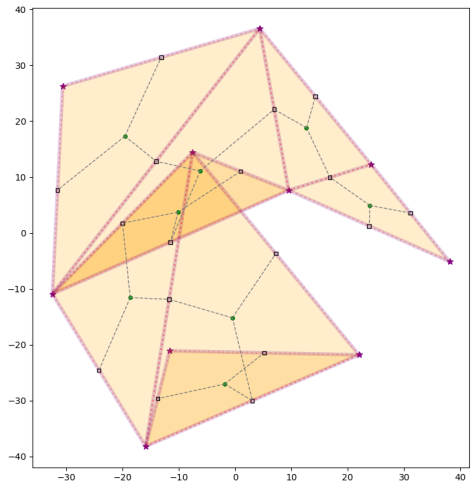

    \centering
    \safeincludegraphics[width=1.0\linewidth]{random\_symmetric\_output.png}
    \caption{Random non-metric instance ($N=10$, visualized by MDS). The formulation relies only on topological connectivity, allowing it to solve instances where the triangle inequality is violated.}
    \label{fig:random_nonmetric}
\end{figure}

To test the topological generality of our formulation, we generated random symmetric non-metric instances of small size (including $N=10$) where edge weights do not satisfy the triangle inequality.
The input was $\mathcal{K}_{All}$. The solver successfully extracted valid dual trees and primal tours; Fig.~\ref{fig:random_nonmetric} shows a representative $N=10$ example. We verified optimality against a standard TSP formulation (MTZ).

\subsection{Comparison}
We compare our surface formulation against the standard Lifted-MTZ baseline on identical sparse inputs (see Appendix \ref{apx:detailed_results} for detailed setup and full results). The surface-based model typically reduces branch-and-bound (B\&B) effort and can solve larger instances more reliably within the class of compact, static MILPs considered. Among surface-based variants, enforcing the global tree via single-commodity flow yields the tightest relaxation and dominates the MTZ-style tree enforcement as $N$ increases.
On TSPLIB Euclidean instances where both baselines are available (\texttt{berlin52}, \texttt{st70}, \texttt{ch130}) under Delaunay sparsification, the Surface-Based-Flow model solves at the root (1 B\&B node), while the Lifted-MTZ baseline requires 2{,}936--16{,}058 nodes (Appendix~\ref{apx:detailed_results}).

\subsection{When does the formulation work well?}
\label{sec:when_works}
The formulation is exact when the candidate triangle set $\mathcal{T}$ is the complete complex $\binom{V}{3}$. In practice, performance heavily relies on restricting $\mathcal{T}$ to a sparse candidate complex.

Empirically, the approach works best when $\mathcal{T}$ is a planar triangulation that contains at least one Hamiltonian tour.
Delaunay triangulations often provide a strong default in Euclidean instances, while other geometric candidates (e.g., Greedy triangulations) can improve solution quality depending on the problem type.

\subsection{Degenerate Inputs}\label{sec:degenerate_inputs}
We further tested “hard-to-solve” instances~\cite{hougardy2018} and degenerate inputs~\cite{hougardy2014} (Appendix \ref{apx:hard_instances}).
\noindent Figure~\ref{fig:hard_instances} visualizes representative solutions for these instances.
On \texttt{Tnm199}, our formulation effectively exploits the richer Greedy triangulation to achieve a 1.85\% gap (vs.\ 6.17\% for Delaunay), surpassing the 5.69\% gap of the Christofides heuristic \cite{christofides1976}.
This confirms that the topological approach can identify superior tours when provided with higher-quality geometric candidates. On \texttt{p100} parallel-line degeneracies, triangulation-based methods degrade badly (14\%--119\% gap vs.\ 2.52\%), exposing a limitation of sparse candidate sets on near-degenerate geometry.

\section{Conclusion}
We have presented a surface-based formulation for TSP that shifts the decision space from edges to triangles.
This changes the optimization object from a primal cycle to a dual surface: instead of directly enforcing Hamiltonicity on $G$, we enforce disk-admissibility of a 2-complex whose boundary is the tour.
Our experiments indicate that this topological shift can yield improved B\&B behavior relative to compact, static edge-based baselines on sparse geometric inputs, particularly when enforcing the global tree constraint via flow.
While exact convergence becomes computationally challenging for larger instances in our current setup, the results suggest that the surface-based model provides a useful foundation for integrating additional geometric or topological constraints that are cumbersome in primal edge formulations, without changing the core model.

\section*{Acknowledgements}
We thank Prof. Stefan Hougardy for his feedback and for providing the test instances. We also thank Anh Quyen Vuong for his suggestions regarding the formulation structure. Finally, we extend our special thanks to the developers of HiGHS for their commitment to providing a powerful and open-source optimization technology.

\appendix
\setcounter{secnumdepth}{3}
\renewcommand{\thesection}{\Alph{section}}

\section{Detailed Computational Results (Standard Instances)}
\label{apx:detailed_results}
\phantomsection

\subsection{Detailed Setup and Performance Analysis}
To evaluate both the topological properties and the solver-specific performance, we benchmarked three formulations:
\begin{enumerate}
    \item \textbf{Conventional:} Standard Lifted-MTZ formulation~\cite{desrochers1991} restricted to edges in the Delaunay triangulation.
    \item \textbf{Surface-Based (MTZ):} Topological formulation using MTZ-based tree constraints~\cite{miller1960}.
    \item \textbf{Surface-Based (Flow):} Topological formulation using Single-Commodity Flow tree constraints~\cite{gavish1978traveling}.
\end{enumerate}

\textbf{Model Selection:} For experiments using sparse planar inputs (Delaunay and Greedy), we utilized the \textbf{Planar Simplification} model described in Appendix \ref{apx:planar}.
This reduces the problem size by eliminating explicit edge node variables ($y_e$), as the incidence graph collapses into the standard dual graph.
We utilized \textbf{HiGHS v1.12}~\cite{huangfu2018} on an \textbf{Apple M4 Max}.
Unless noted, runtime results are arithmetic means of 10 runs with distinct random seeds.
Problem sizes (Vars/Cons) are reported after solver presolve.
Table \ref{tab:combined_results} summarizes the performance across all instances.
The results reveal a clear transition in which the \textbf{Surface-Based-Flow} formulation dominates \textbf{Surface-Based-MTZ} as $N$ grows.
For small instances ($N \le 130$), both surface-based models typically solve at the root.
For larger instances, the single-commodity flow constraint yields a substantially tighter relaxation, reflected in dramatically smaller Branch-and-Bound trees (e.g., 142 nodes vs.\ 23{,}205 for \texttt{pcb442}).
This supports the interpretation that our formulation benefits from flow-based connectivity enforcement more strongly than from MTZ-style ordering constraints, even though both encode the same topological admissibility conditions.
We compare against compact static models, excluding dynamic cut-generation methods (e.g., Concorde~\cite{applegate2006}) which remain state-of-the-art for pure TSP.

\begin{table*}[htbp]
    \centering
    \caption{\textbf{Combined Computational Results.} Post-presolve MILP size and performance; values are 10-run means (unless noted), with a timeout limit of 1-hour.}
    \label{tab:combined_results}
    \footnotesize
    \setlength{\tabcolsep}{5pt}
    \renewcommand{\arraystretch}{1.0}
    \begin{tabular}{l l | rr | rr}
        \hline
        & & \multicolumn{2}{c|}{\textbf{MILP Size}} & \multicolumn{2}{c}{\textbf{HiGHS}} \\
        \textbf{Instance} & \textbf{Model} & \textbf{Vars} & \textbf{Cons} & \textbf{Time (s) [Solved]} & \textbf{Nodes} \\
        \hline
        \textbf{Berlin52} & Lifted-MTZ & 341 & 523 & 0.08s [10/10] & 1 \\
        ($N=52$) & Surface-Based-MTZ & 553 & 1,107 & \textbf{0.02s} [10/10] & 1 \\
         & Surface-Based-Flow & 684 & 1,054 & 0.03s [10/10] & \textbf{1} \\
        \hline
        \textbf{St70} & Lifted-MTZ & 463 & 723 & 14.23s [10/10] & 2,936 \\
        ($N=70$) & Surface-Based-MTZ & 759 & 1,521 & 4.35s [10/10] & 257 \\
         & Surface-Based-Flow & 926 & 1,440 & \textbf{0.20s} [10/10] & \textbf{1} \\
        \hline
        \textbf{Ch130} & Lifted-MTZ & 879 & 1,362 & 108.6s [10/10] & 16,058 \\
        ($N=130$) & Surface-Based-MTZ & 1,458 & 2,922 & \textbf{0.29s} [10/10] & \textbf{1} \\
         & Surface-Based-Flow & 1,793 & 2,775 & 1.62s [10/10] & \textbf{1} \\
        \hline
        \textbf{Tnm199} & Surface-Based-Flow (Del) & 2,764 & 4,283 & \textbf{0.51s} [10/10] & \textbf{1} \\
        ($N=199$) & Surface-Based-Flow (Grd) & 2,090 & 3,202 & 4.76s [10/10] & 1.6 \\
        \hline
        \textbf{A280} & Lifted-MTZ & 1,907 & 2,992 & DNF [0/5] & $>$53k \\
        ($N=280$) & Surface-Based-MTZ & 3,188 & 6,379 & 29.0s [10/10] & 490 \\
         & Surface-Based-Flow & 3,923 & 6,056 & \textbf{5.68s} [10/10] & \textbf{1} \\
        \hline
        \textbf{Pcb442} & Surface-Based-MTZ & 5,195 & 10,393 & 1656s* [7/10] & 23,205 \\
        ($N=442$) & Surface-Based-Flow & 6,412 & 9,879 & \textbf{81.1s} [10/10] & \textbf{142} \\
        \hline
        \textbf{d657} & Surface-Based-Flow & 9,634 & 14,837 & TLE (0.60\%) & $>$30k \\
        ($N=657$) & & & & & \\
        \hline
    \end{tabular}
    \par\vspace{0.1cm}
    {\footnotesize * For \texttt{pcb442} Surface-Based-MTZ, 3/10 runs exceeded the 1-hour time limit; time and nodes are averaged over the 7 solved runs.}
\end{table*}

\subsection{Impact of Input Sparsity on Lifted-MTZ}
\label{apx:sparsity}
To justify the choice of Delaunay-restricted input sets for the baseline comparison, we evaluated the performance of the standard Lifted-MTZ formulation on the complete graph $K_N$ (Table \ref{tab:complete_graph}).
The results show a massive divergence in performance.
On \texttt{st70}, HiGHS runtime increases by roughly 13.5x (14.2s vs.\ 191.7s) when moving from the Delaunay-restricted graph to the complete graph $K_N$.
The topological formulation remains extremely fast on the sparse complex (0.2s), indicating that the topological structure and the restricted geometry both contribute materially to scalability (a further $\sim$71x speedup).

\begin{table}[htbp]
    \centering
    \caption{\textbf{Lifted-MTZ on Complete Graph vs. Delaunay.} On \texttt{st70} ($N=70$), sparsity gives $\sim$13.5x and the surface-based approach gives a further $\sim$71x.}
    \label{tab:complete_graph}
    \small
    \setlength{\tabcolsep}{5pt}
    \renewcommand{\arraystretch}{1.1}
    \begin{tabular}{l l r}
        \hline
        \textbf{Instance} & \textbf{Edge Set} & \textbf{HiGHS (s)} \\
        \hline
        \textbf{St70} & Complete ($K_N$) & 191.7s \\
         & Delaunay & 14.2s \\
         & Topological & \textbf{0.2s} \\
        \hline
    \end{tabular}
\end{table}

\section{Extended Experimental Results (Hard Instances)}
\label{apx:hard_instances}
\phantomsection

\begin{figure}[b!]
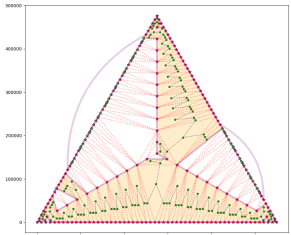
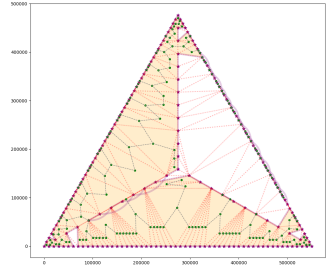
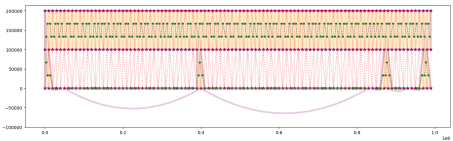
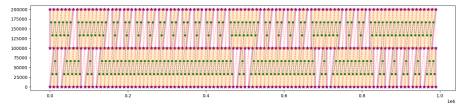

    \centering
    \begin{subfigure}[b]{0.49\textwidth}
        \centering
        \safeincludegraphics[width=\linewidth]{Tnm199\_delaunay\_output.png}
        \caption{\texttt{Tnm199} (Delaunay): Gap 6.17\%}
    \end{subfigure}
    \hfill
    \begin{subfigure}[b]{0.49\textwidth}
        \centering
        \safeincludegraphics[width=\linewidth]{Tnm199\_greedy\_output.png}
        \caption{\texttt{Tnm199} (Greedy): Gap 1.85\%}
    \end{subfigure}

    \vspace{1em}

    \begin{subfigure}[b]{0.49\textwidth}
        \centering
        \safeincludegraphics[width=\linewidth]{p100\_delaunay\_output.png}
        \caption{\texttt{p100} (Delaunay): Gap 14.5\%}
    \end{subfigure}
    \hfill
    \begin{subfigure}[b]{0.49\textwidth}
        \centering
        \safeincludegraphics[width=\linewidth]{p100\_greedy\_output.png}
        \caption{\texttt{p100} (Greedy): Gap 119.6\%}
    \end{subfigure}

    \caption{Solution quality on Hard Instances. Top Row: On \texttt{Tnm199}, Greedy (b) finds a superior surface to Delaunay (a). Bottom Row: On the degenerate \texttt{p100}, heuristics fail; Greedy (d) is especially catastrophic compared to Delaunay (c).}
    \label{fig:hard_instances}
\end{figure}

To further evaluate the impact of the input complex on solution quality, we tested the formulation on particularly difficult instances.

\subsection{Analysis of Hard Instances: Tnm199}
We tested on \texttt{Tnm199} ($N=199$) from the "Hard-to-Solve" dataset~\cite{hougardy2018}. This instance is constructed to be practically intractable for standard exact solvers; Concorde requires 411,222 seconds ($\approx 4.7$ days) to find the global optimum. We compared our topological method using two input triangulation strategies (Delaunay and Greedy) against the standard Christofides algorithm. The results (Table \ref{tab:tnm199_results}) illustrate the model's dependence on input quality. While the Delaunay-based solution (Fig~\ref{fig:hard_instances}a) is structurally limited (similar to Christofides), the Surface-Based-Flow model on the Greedy triangulation successfully extracts the superior geometry within the complex (Fig.~\ref{fig:hard_instances}b), narrowing the gap to \textbf{1.85\%}. This demonstrates that the topological formulation scales in quality with the richness of the input complex, unlike fixed constructive heuristics.

\begin{table}[htbp]
    \centering
    \caption{\textbf{Results on Tnm199.} Greedy input enables a tighter gap than Delaunay and Christofides.}
    \label{tab:tnm199_results}
    \small
    \setlength{\tabcolsep}{4pt}
    \renewcommand{\arraystretch}{1.1}
    \begin{tabular}{l r r r}
        \hline
        \textbf{Method} & \textbf{Time (s)} & \textbf{Objective} & \textbf{Gap (\%)} \\
        \hline
        Global Opt. (Concorde) & $>$400k & 3,139,778 & 0.00 \\
        \hline
        Christofides & $<$ 1.00 & 3,318,561 & +5.69 \\
        Surface-Based-Flow (Delaunay) & 0.51 & 3,333,452 & +6.17 \\
        Surface-Based-Flow (Greedy) & 4.76 & 3,197,767 & \textbf{+1.85} \\
        \hline
    \end{tabular}
\end{table}

\subsection{Sensitivity Analysis: Parallel Line Degeneracy}
To test the limitations of sparsification under degenerate geometry, we evaluated the instance \texttt{p100.100000} ($N=100$), which arranges points on three parallel lines~\cite{hougardy2014}. We compared the Surface-Based-Flow formulation (using Delaunay and Greedy inputs) against the standard Christofides heuristic. The results in Table \ref{tab:p100_results} reveal a critical vulnerability in geometric sparsifiers. While Christofides remains robust (2.52\% average gap), the triangulation-based methods suffer substantially (Fig.~\ref{fig:hard_instances}c-d).

\begin{table}[!ht]
    \centering
    \caption{\textbf{Results on p100.100000.} Parallel-line degeneracy breaks triangulation-based candidate sets.}
    \label{tab:p100_results}
    \small
    \setlength{\tabcolsep}{4pt}
    \renewcommand{\arraystretch}{1.1}
    \begin{tabular}{l r r}
        \hline
        \textbf{Method} & \textbf{Objective} & \textbf{Gap (\%)} \\
        \hline
        \textbf{Global Opt.} (Exact) & 4,160,200 & 0.00 \\
        \hline
        Christofides & 4,265,038 & \textbf{+2.52} \\
        Surface-Based-Flow (Delaunay) & 4,762,994 & +14.49 \\
        Surface-Based-Flow (Greedy) & 9,136,966 & +119.63 \\
        \hline
    \end{tabular}
\end{table}

\section{Explicit MILP Formulation}
\label{apx:full_model}
\phantomsection
Here we provide the complete Mixed-Integer Linear Program corresponding to the high-level formulation in Section \ref{sec:formulation}.

\subsection{Logical Linking and Closure}
Active simplices must be connected via the triangle--edge incidence relation $\mathcal{A}$.
\begin{align}
    & z_{te} \le x_t, \quad z_{te} \le y_e \quad \forall (t,e) \in \mathcal{A}. \label{eq:linking}
\end{align}

Define the induced degrees
\[
    \deg_U(t) := \sum_{e:(t,e)\in\mathcal{A}} z_{te} \quad \forall t\in U,
    \qquad
    \deg_W(e) := \sum_{t:(t,e)\in\mathcal{A}} z_{te} \quad \forall e\in W.
\]

Degree/Closure Constraint: If a triangle is active, all three of its edges must be active:
\begin{equation}
\deg_U(t) = 3x_t \quad \forall t \in U.
\end{equation}

\subsection{Topological Regularity}
These global counts follow from the Euler characteristic of the intended disk-like complex~\cite{edelsbrunner2010}.
\begin{align}
    & \sum_{t \in U} x_t = N - 2 \quad (\text{Triangle Cardinality}) \label{eq:card_tri} \\
    & \sum_{e \in W} y_e = 2N - 3 \quad (\text{Edge Cardinality}) \label{eq:card_edge} \\
    & y_e \le \deg_W(e) \le 2y_e \quad \forall e \in W \quad (\text{Degree/Manifold}). \label{eq:manifold_apx}
\end{align}

\subsection{Local Euler Filter and Connectivity}
For each city $v \in V$, let
\[
U(v) := \{t \in U \mid v \in t\},
\qquad
W(v) := \{e \in W \mid v \in e\},
\]
and define the restricted incidence set
\[
\mathrm{Conn}(v) := \{(t,e) \in \mathcal{A} \mid t \in U(v),\ e \in W(v)\}.
\]
Define the local counts
\[
T(v) := \sum_{t\in U(v)} x_t,
\]
\[
E(v) := \sum_{e\in W(v)} y_e,
\]
\[
I(v) := \sum_{(t,e)\in \mathrm{Conn}(v)} z_{te},
\]
and the local Euler characteristic
\[
\chi(v) := T(v) + E(v) - I(v).
\]
We enforce~\cite{edelsbrunner2010}
\begin{equation}
    \chi(v) = 1 \quad \forall v \in V.
    \label{eq:euler_milp}
\end{equation}

\subsection{Global Tree Constraint (Directed Single-Commodity Flow)}
We designate a root node $r$ dynamically. To break symmetries and improve performance, we restrict the candidate root set $R_{cand} \subset U$ to the set of triangles incident to the city with the minimum
degree in the mesh. We introduce binary variables $r_t \in \{0,1\}$ for $t \in R_{cand}$ and enforce:
\begin{equation}
    \sum_{t \in R_{cand}} r_t = 1,
    \qquad
    r_t \le x_t \ \ \forall t\in R_{cand}.
\end{equation}
For all other nodes $i \notin R_{cand}$ (including all edge nodes $W$), we define $r_i := 0$. To enforce linearity and tighter relaxations, we model a
\textit{Rooted Directed Tree} on the bipartite graph. Let
\[
A_{\text{dir}} = \{(t,e),(e,t) \mid (t,e) \in \mathcal{A}\}
\]
be the set of directed arcs. The total number of active nodes is a known constant
\[
K = (N-2) + (2N-3) = 3N-5.
\]

\paragraph{Active-node indicator}
For each node $i \in U \cup W$, define
\[
\text{active}_i :=
\begin{cases}
x_t & \text{if } i=t \in U,\\
y_e & \text{if } i=e \in W.
\end{cases}
\]

\textbf{1. Directed Arc Selection:}
For every incidence $(t,e)\in\mathcal{A}$, we introduce two directed binary variables $a_{te}, a_{et} \in \{0,1\}$ and enforce:
\begin{equation}
    a_{te} + a_{et} = z_{te} \quad \forall (t,e) \in \mathcal{A}.
\end{equation}

\textbf{2. Flow Capacity and Conservation:}
Let $f_{uv} \ge 0$ be the continuous flow on arc $(u,v)\in A_{\text{dir}}$. Let $M$ be a sufficiently large constant (e.g., $M=K$).
\begin{align}
    & f_{uv} \le M \cdot a_{uv}
    \quad \forall (u,v) \in A_{\text{dir}}
    \quad (\text{Capacity}) \label{eq:flow_cap} \\
    & \sum_{u:(u,i)\in A_{\text{dir}}} f_{ui} - \sum_{w:(i,w)\in A_{\text{dir}}} f_{iw}
      = \text{active}_i - K \cdot r_i
      \quad \forall i \in U \cup W
    \quad (\text{Conservation}) \label{eq:flow}
\end{align}

\textbf{3. Rooted Tree Structure:}

Every active non-root node has exactly one incoming arc:
\begin{equation}
    \sum_{u: (u,i) \in A_{\text{dir}}} a_{ui} = \text{active}_i - r_i
    \quad \forall i \in U \cup W.
    \label{eq:in_arcs}
\end{equation}

\section{Special Case: The Planar Simplification}
\label{apx:planar}
\phantomsection

\begin{figure}[b!]
\centering
\begin{tikzpicture}[scale=1.0, >=Stealth]
    \tikzset{arrow_edge/.style={decoration={markings, mark=at position 0.6 with {\arrow{>}}}, postaction={decorate}},
             cancel_edge/.style={draw=red, thick, dashed},
             boundary_edge/.style={draw=magenta, line width=2pt},
             surface/.style={fill=orange, fill opacity=0.3}}
    \begin{scope}[shift={(0,0)}]
        \coordinate (L1) at (0,0); \coordinate (L2) at (2,0); \coordinate (L3) at (1,1.732);
        \fill[surface] (L1) -- (L2) -- (L3) -- cycle;
        \draw[arrow_edge] (L1) -- (L2); \draw[arrow_edge] (L2) -- (L3); \draw[arrow_edge] (L3) -- (L1);
        \node at (1, 0.6) {$t_1$};
        \begin{scope}[shift={(2.5,0)}]
            \coordinate (R1) at (0,0); \coordinate (R2) at (-1,1.732);
            \coordinate (R3) at (1,1.732);
            \fill[surface] (R1) -- (R2) -- (R3) -- cycle;
            \draw[arrow_edge] (R1) -- (R2); \draw[arrow_edge] (R2) -- (R3);
            \draw[arrow_edge] (R3) -- (R1);
            \node at (0, 0.6) {$t_2$};
        \end{scope}
    \end{scope}
    \draw[->, ultra thick, gray] (4.5, 1) -- (5.5, 1) node[midway, above] {Glue};
    \begin{scope}[shift={(6,0)}]
        \coordinate (V1) at (0,0);
        \coordinate (V2) at (2,0); \coordinate (V3) at (1,1.732); \coordinate (V4) at (3, 1.732);
        \fill[surface] (V1) -- (V2) -- (V3) -- cycle; \fill[surface] (V2) -- (V4) -- (V3) -- cycle;
        \draw[boundary_edge] (V1) -- (V2);
        \draw[boundary_edge] (V2) -- (V4);
        \draw[boundary_edge] (V4) -- (V3); \draw[boundary_edge] (V3) -- (V1);
        \draw[cancel_edge] (V2) -- (V3);
        \draw[->, red, thick] ($(V2)!0.3!(V3)$) -- ($(V2)!0.7!(V3)$) node[midway, left, font=\tiny] {$e_{1}$};
        \draw[->, red, thick] ($(V3)!0.3!(V2)$) -- ($(V3)!0.7!(V2)$) node[midway, right, font=\tiny] {$e_{2}$};
        \node at (1.5, -0.8) {$- 2 \times \text{Length}(\text{Internal})$};
        \coordinate (C_T1) at (barycentric cs:V1=1,V2=1,V3=1); \coordinate (C_T2) at (barycentric cs:V2=1,V3=1,V4=1);
        \draw[blue!80!black, thick, dashed] (C_T1) -- (C_T2);
        \fill[blue!80!black] (C_T1) circle (1.5pt); \fill[blue!80!black] (C_T2) circle (1.5pt);
    \end{scope}
\end{tikzpicture}
\caption{Combinatorial boundary identity: internal shared edges cancel by opposite orientations, leaving only the boundary.}
\label{fig:stokes_cancellation}
\end{figure}
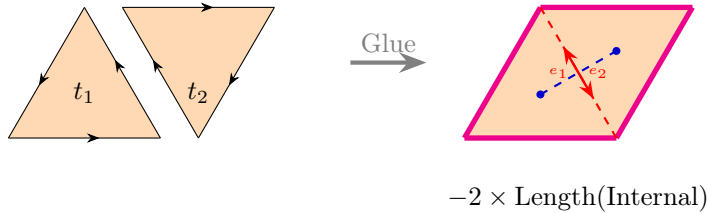

In the general case (Section \ref{sec:formulation}), the Bipartite Graph is necessary because a primal edge $e$ (edge node) may be incident to many triangles. However, if the candidate set $\mathcal{T}$ is a planar triangulation (e.g., Delaunay), every primal edge is incident to at most 2 triangles.
In this specific case, the explicit edge nodes $y_e$ can be removed. The Bipartite Graph $B$ collapses into a simple dual graph, where edges connect neighboring triangles (Fig.~\ref{fig:stokes_cancellation}).
Consequently, the Euler filter simplifies from the bipartite form ($\chi(v)=1$) to the standard dual-graph form ($|V_{dual}| - |E_{dual}| = 1$).
This reduces the variable count by $O(|E|)$ and simplifies the constraints. The manifold degree constraint becomes implicit.

\paragraph{Boundary Cancellation}
In this planar case, the objective simplifies by eliminating explicit edge-node variables, and can be derived from Eq.~(\ref{eq:objective_main}).
\begin{equation}
    \text{Min } Z = \sum_{t \in \mathcal{T}_{selected}} P(t) - 2 \sum_{e \in \mathcal{E}_{dual}} L(e)
\end{equation}
where $P(t)$ is the perimeter of triangle $t$ and $L(e)$ is the length of the primal edge shared by the two triangles connected by dual edge $e$.

\end{document}